\pdfoutput=1
\RequirePackage{ifpdf}
\ifpdf 
\documentclass[pdftex]{sigma}
\else
\documentclass{sigma}
\fi

\numberwithin{equation}{section}

\newtheorem{Theorem}{Theorem}
\newtheorem*{Theorem*}{Theorem}
\newtheorem{Corollary}[Theorem]{Corollary}
\newtheorem{Lemma}[Theorem]{Lemma}

\theoremstyle{definition}

\newtheorem{Example}[Theorem]{Example}

\usepackage{tikz-cd}

\begin{document}

\allowdisplaybreaks

\newcommand{\arXivNumber}{2601.07714}

\renewcommand{\thefootnote}{}

\renewcommand{\PaperNumber}{051}

\FirstPageHeading

\ShortArticleName{A Note on Somewhere Positive Loops of Contactomorphisms}

\ArticleName{A Note on Somewhere Positive Loops\\ of Contactomorphisms\footnote{This paper is a~contribution to the Special Issue on Geometry and Dynamics in memory of Will Merry. The~full collection is available at \href{https://sigma-journal.com/Merry.html}{https://sigma-journal.com/Merry.html}}}

\Author{Igor ULJAREVI\'C}

\AuthorNameForHeading{I.~Uljarevi\'c}

\Address{Faculty of Mathematics, University of Belgrade,\\ Studentski trg 16, 11158 Belgrade, Republic of Serbia}
\Email{\mail{igor.uljarevic@matf.bg.ac.rs}}

\ArticleDates{Received January 13, 2026, in final form May 06, 2026; Published online May 18, 2026}

\Abstract{In this note, we consider contractible loops of contactomorphisms that are positive over some non-empty closed subset of a contact manifold. Such closed subsets are called immaterial. We argue that the complement of a Reeb-invariant immaterial subset can be seen as big in contact geometric terms. This is supported by two results: one regarding symplectic homology of the filling and the other regarding recently introduced contact quasi-measures.}


\Classification{53D10}

\renewcommand{\thefootnote}{\arabic{footnote}}
\setcounter{footnote}{0}

\section{Introduction}

A subset $A$ of a cooriented contact manifold $M$ is called \emph{immaterial}\footnote{The notion was introduced in \cite[Section~7]{uljarevic2023selective}.} if there exists a contractible loop of contactomorphisms $\varphi_t\colon M\to M$ whose vector field $X$ points in the positive direction (of the contact distribution) on $A$.
The property of being immaterial is clearly related to orderability, a well studied phenomenon in contact geometry originally introduced in \cite{eliashberg2000partially}. Namely, a contact manifold~$M$ is non-orderable if, and only if, all of its subsets are immaterial. The existence of non-orderable contact manifolds shows that immaterial subsets are not necessarily small in the sense of contact squeezing or contact displaceability. Indeed, the whole contact manifold~$M$ (within itself) cannot be squeezed into a smaller subset or displaced from itself, although the whole $M$ is immaterial if non-orderable.
The only immaterial subset in dimension 1 is the empty set.\footnote{{Here, we prove the claim for $\mathbb{S}^1$. It is enough to show that a point $p\in\mathbb{S}^1$ is not immaterial. Let $\varphi_t\colon \mathbb{S}^1\to \mathbb{S}^1$ be a loop of contactomorphisms that is positive at $p$. Then, the degree of the map $\gamma\colon\mathbb{S}^1\to \mathbb{S}^1$, $t\mapsto \varphi_t(p)$ is positive. (This can be seen by counting elements of $\gamma^{-1}(\{p\})$ with signs.) Therefore, $\{\varphi_t\}_t$ cannot be contractible.}} This is in contrast with dimensions greater than 3 where every compact subset of a Darboux chart is immaterial \cite[Lemma~7.3]{uljarevic2023selective}.

{The results of this note can be seen as contributions to questions regarding the size in contact geometry. These questions are quite delicate (as evident from work by Eliashberg--Kim--Polterovich \cite{eliashberg2006geometry}) and have been an active research area over the last two decades. Our first result (Theorem~\ref{thm:complement-big}) is formulated in terms of selective symplectic homology \cite{uljarevic2023selective}, while the second (Theorem~\ref{thm:tau}) is expressed using recently introduced contact quasi-measures \cite{djordjevic2025quantitative}. The contact quasi-measures from \cite{djordjevic2025quantitative} are strongly inspired by Entov--Polterovich's construction of (symplectic) quasi-measures \cite{EP06} using spectral invariants from \cite{oh2004construction}. See \cite{Ent14} or \cite{polterovich2014function} for a more detailed historical account.}

{A central role in this paper is played by \emph{contact Hamiltonian Floer homology}, introduced in~\cite{merry2019maximum} by proving a maximum principle for slopes that are allowed to be general contact Hamiltonians (instead of constants in $\mathbb{R}$ as common in symplectic homology theory). The cases where the contact Hamiltonian (i.e., the slope) is positive \cite{seidel2008biased} or Reeb-invariant \cite{ritter2016circle} were known prior to \cite{merry2019maximum}. The theory of contact Hamiltonian Floer homology was further developed in \cite{djordjevic2023quantitative, uljarevic2022hamiltonian} and in \cite{cant2023shelukhin, cant2024remarks} with independent contributions.}

\subsection{Immaterial subsets and symplectic homology}

In this section, $W$ denotes a Liouville domain.\footnote{We expect the result of this section to hold in a more general framework than that of Liouville domains. See \cite[Remark~1.5]{uljarevic2023selective}.} Selective symplectic homology $SH_\ast^\Omega(W)$ is a~version of symplectic homology introduced in \cite{uljarevic2023selective}. It is obtained as the direct limit of Floer homology groups for Hamiltonians whose slopes tend to $+\infty$ on the conical end above the open subset $\Omega\subset \partial W $ of the boundary, but remain close to 0 and positive elsewhere. The more precise definition of $SH_\ast^\Omega(W)$ is recalled in Section~\ref{sec:proofSH}.
If $\Omega=\partial W$, then $SH^\Omega_\ast(W)$ coincides with Viterbo's {\cite{viterbo1999functors}} symplectic homology $SH_\ast(W)$. Given two open subsets $\Omega_1\subset \Omega_2$ of $\partial W$, there is a functorial morphism $SH^{\Omega_1}_\ast(W)\to SH^{\Omega_2}_\ast(W)$ called \emph{continuation map}. In particular, the continuation map $SH_\ast^\Omega(W)\to SH_\ast(W)$ is well defined for every open subset $\Omega\subset\partial W$.
The following theorem shows that complements of immaterial Reeb-invariant subsets are big in terms of selective symplectic homology.

\begin{Theorem}\label{thm:complement-big}
Let $A\subset \partial W$ be the preimage of a closed subset under a smooth Reeb-invariant map {$\partial W\to\mathbb{R}$}. Assume $A$ is immaterial and denote $\Omega:=\partial W\setminus A$. Then, the continuation map
\begin{equation*}SH^{\Omega}_\ast(W)\to SH_\ast(W)\end{equation*}
is surjective.
\end{Theorem}

Theorem~\ref{thm:complement-big} can be seen as a generalization and an improvement of \cite[Theorem~7.4]{uljarevic2023selective}. {Namely, \cite[Theorem~7.4]{uljarevic2023selective} claims that the rank of the map $SH^{\Omega}_\ast(W)\to SH_\ast(W)$ is equal to $\dim SH_\ast(W)$ in the case where $A$ is a transversely embedded circle. It is an interesting question whether the condition of $A$ being the preimage of a closed subset under a smooth Reeb-invariant map~${\partial W\to \mathbb{R}}$ can be relaxed.}

\subsection{Contact quasi-measures and immaterial subsets}

In this section, $W$ stands for a weakly\textsuperscript{+} monotone strong filling of a closed contact manifold such that the unit $e\in SH_\ast(W)$ is not eternal. {In what follows, we recall the notions of a weakly\textsuperscript{+} monotone symplectic manifold and of an eternal element. A} $2n$-dimensional symplectic manifold $(W, \omega)$ is weakly\textsuperscript{+} monotone\footnote{ See \cite[Lemma~1.1]{hofer1995floer}, \cite[Assumption ($W^+$)]{seidel1997pi_1}, and \cite[Section~2B]{ritter2016circle}.} if at least one of the following conditions is satisfied:
\begin{enumerate}\itemsep=0pt
\item[(1)] $ \omega |_{\pi_2(W)}=0$,
\item[(2)] $ c_1 |_{\pi_2(W)}=0$,
\item[(3)] there exists $s>0$ such that\footnote{Here, both $\omega$ and $c_1$ are seen as maps assigning a number to an element of $\pi_2(W)$.} $ \omega |_{\pi_2(W)}= s\cdot c_1 |_{\pi_2(W)}$,
\item[(4)] the minimal Chern number of $W$ is at least $n-1$.
\end{enumerate}
These conditions ensure that the symplectic homology of $W$ is well defined. Now, we also recall the notion of being eternal \cite{cant2024remarks}. An element $\theta\in SH_\ast(W)$ is called {\emph{eternal}} if it is contained in the image of the canonical morphism
\[ HF_\ast(H^a)\to SH_\ast(W) \]
for all admissible\footnote{{Recall that a number $a\in\mathbb{R}$ is \emph{admissible} if it is not the period of a closed Reeb orbit on $\partial W$. The extension of this notion to contact Hamiltonians is recalled on page~\pageref{def:admissible}}.} $a\in \mathbb{R}$. Here \smash{$H^a_t\colon \widehat{W}\to \mathbb{R}$} denotes a non-degenerate Hamiltonian on the completion \smash{$\widehat{W}$} of $W$ having slope $a$ at the conical end.

Inspired by work of Entov and Polterovich \cite{Ent14, EP06,EP08,EP-rigidity}, the paper \cite{djordjevic2025quantitative} introduced {a map~$\tau$, called \emph{partial contact quasi-measure},} which assigns a number from $[0,1]$ to a closed subset of~$\partial W$. The quasi-measure $\tau$ is constructed using contact spectral invariants $c(h, \theta)$ introduced independently and simultaneously in \cite{cant2023shelukhin} and \cite{djordjevic2023quantitative}. We recall the construction of $\tau$ in Section~\ref{sec:tau}. The following theorem asserts that each contractible loop of contactomorphisms furnishes a lower bound for~$\tau$.

\begin{Theorem}\label{thm:tau}
Let $W$ and $\tau$ be as above, and $B\subset \partial W$ a closed subset. Let $f\colon\partial W\times [0,1]\to\mathbb{R}$ be a contact Hamiltonian giving rise to a contractible loop of contactomorphisms. Denote {by}~$m$ and $m_B$ the infima of $f$ over $\partial W\times[0,1]$ and $(\partial W\setminus B)\times[0,1]$, respectively. Then, the following inequality holds:
\[\tau(B)\geqslant \frac{m_B}{m_B-m}\]
if $m_B\not=m$.
\end{Theorem}

{By definition, $\tau(B)$ is an element of $[0,1]$. In this paragraph, we explain why Theorem~\ref{thm:tau} does not contradict this fact. The number $m$ cannot be positive, because that would imply that $\partial W$ is not orderable and contradict \cite[Theorem~1.6]{djordjevic2025quantitative}. Therefore, $m_B-m\geqslant m_B$ and $\frac{m_B}{m_B-m}\leqslant 1$. This means that Theorem~\ref{thm:tau} never forces $\tau(B)$ to be greater than 1.

It is useful to observe that Theorem~\ref{thm:tau} is only interesting in the case where $\partial W\setminus B$ is immaterial. Indeed, if $m_B$ is negative, then the theorem asserts that $\tau(B)$ is greater than some negative number, and that follows directly from the definition of $\tau$.}
Lemma~9.3 from \cite{djordjevic2025quantitative} shows $\tau(B)=0$ if $B$ is displaceable and the preimage of a closed subset under a Reeb-invariant function. Therefore, Theorem~\ref{thm:tau} and \cite[Lemma~9.3]{djordjevic2025quantitative} imply the following corollary.

\begin{Corollary}\label{cor:nondisp}
Let $B\subset\partial W$ be the preimage of a closed subset under a Reeb-invariant smooth function {$\partial W\to\mathbb{R}$}. If $\overline{\partial W \setminus B}$ is immaterial, then $B$ is contact non-displaceable.
\end{Corollary}

It is an interesting question whether the condition of $B$ being the preimage of a closed subset under a smooth Reeb-invariant map $\partial W\to \mathbb{R}$ can be relaxed.
The next example applies Theorem~\ref{thm:tau} in a concrete situation of $S^\ast \mathbb{S}^n$. The example illustrates a general method of finding immaterial subsets in unit cotangent bundles by lifting contractible loops of diffeomorphisms from the base.

\begin{Example}
Let $n\in\mathbb{N}$ be a number equal to 3 modulo 4 and let $\psi_t\colon \mathbb{S}^n\to\mathbb{S}^n$ be given by $\psi_t(x):= {\rm e}^{2\pi {\rm i} t}\cdot x$ (we see $\mathbb{S}^n$ as the unit sphere in the complex space). Since $n\geqslant 3$ is equal to 3 modulo 4, $\psi_t$ is a contractible loop of diffeomorphisms.\footnote{The loop $\psi_t$ represents an even power of the generator of $\pi_1\bigl(U\bigl(\frac{n+1}{2}\bigr)\bigr)\cong\mathbb{Z}$ because $\frac{n+1}{2}$ is even. Therefore, $\psi_t$ represents the trivial element in $\pi_1(O(n+1))\cong\mathbb{Z}_2$. Consequently, $\psi_t$ is contractible as a loop of diffeomorphisms.} Therefore,
\[
\varphi_t\colon \ S^\ast\mathbb{S}^n\to S^\ast\mathbb{S}^n, \qquad v^\ast\mapsto \psi^\ast_t(v^\ast)
\]
is a contractible loop of contactomorphisms. The contact Hamiltonian of $\varphi_t$ is given by
\[
f\colon \ S^\ast\mathbb{S}^n\to[-2\pi, 2\pi], \qquad f(v^\ast):= v^\ast(X(p)),
\]
where $X$ is the vector field of $\psi_t$ and $p$ is the base point for $v^\ast$. Denote by $B_c$ the subset $\{v^\ast\in S^\ast\mathbb{S}^n\mid v^\ast(X)\leqslant c\}$. The set $B_c$ for $c\in(-2\pi, 2\pi), $ can be seen as a fibration over $\mathbb{S}^n$ with each fibre diffeomorphic to the $(n-1)$-dimensional ball $\mathbb{B}^{n-1}$. Now, Theorem~\ref{thm:tau} implies
\[ \tau(B_c)\geqslant \frac{c}{c+2\pi}. \]
As a consequence of Corollary~\ref{cor:nondisp}, the set $B_c$ is contactly non-displaceable for $c>0$. This claim, however, is not new. In fact, \cite[Example~2.11]{sun2026contact} asserts that $B_c$ is non-displaceable for $c\geqslant 0$.
\end{Example}

\section{Proof of Theorem~\ref{thm:complement-big}}\label{sec:proofSH}

 We first recall the definition of the selective symplectic homology $SH_\ast^\Omega(W)$. For an open subset $\Omega\subset \partial W$, denote by $\mathcal{H}_\Omega$ the set of the smooth autonomous contact Hamiltonians $h\colon \partial W\to[0,+\infty)$ such that\footnote{The second condition is slightly different from that of \cite[Definition~3.1]{uljarevic2023selective}. However, using either of these two conditions yields the same definition of selective symplectic homology.}
 \begin{enumerate}\itemsep=0pt
 \item[(1)] ${\rm supp} h \subset \Omega$,
 \item[(2)] $h$ vanishes up to order $2$ at each $p\in\partial W$ such that $h(p)=0$,
 \item[(3)] the 1-periodic orbits of $h$ are constant.
 \end{enumerate}

 Denote by $\Pi(h)$ the set of the positive smooth functions $f\colon \partial W\to\mathbb{R}_{{>}0}$ such that the contact Hamiltonian $h+f$ has no 1-periodic orbits. {We see $\mathcal{H}_\Omega$ and $\Pi(h)$ as directed sets with the preorder relation pointwise $\leqslant$.} The selective symplectic homology is defined by
 \begin{equation}\label{eq:SSH} SH_\ast^\Omega(W):= \underset{h\in\mathcal{H}_\Omega}{\lim_{\longrightarrow}} \underset{f\in\Pi(h)}{\lim_{\longleftarrow}} HF_\ast(h+f),\end{equation}
 where the limits are taken with respect to the continuation maps {\cite[Section~4]{merry2019maximum}}. Here, $HF_\ast(h)$ stands for the Floer homology for the contact Hamiltonian $h$, that was introduced in \cite{merry2019maximum}. Notice that the inverse limit in \eqref{eq:SSH} stabilizes \cite[Lemma~3.5]{uljarevic2023selective}. In the following lemma and in the rest of the paper, \smash{$\varphi^f_t$} denotes the contact isotopy furnished by the contact Hamiltonian $f\colon\partial W\times\mathbb{R}\to\mathbb{R}$. {The contact Hamiltonian of the composition \smash{$\varphi_t^h\circ\varphi_t^f$} is denoted by $h\# f$. Explicitly,
 \[ (h\# f)_t(x):= h_t(x) + \bigl(\kappa_t^h\cdot f_t\bigr)\circ\bigl(\varphi^h_t\bigr)^{-1}(x), \]
 where $\kappa^h_t$ is defined by \smash{$\bigl(\varphi^h_t\bigr)^\ast\alpha =\kappa_t^h\cdot \alpha$}.}

\begin{Lemma}\label{lem:twisttrick}
In the situation of Theorem~$\ref{thm:complement-big}$, for every $\ell\in\mathbb{R}$ there exists $h\in\mathcal{H}_{\Omega}(\partial W)$ and a~contractible loop \smash{$\varphi^f_t\colon \partial W\to\partial W$} of contactomorphisms\footnote{Here and in the rest of the paper, loops $\varphi_t\colon \partial W\to\partial W$ of contactomorphisms are seen as 1-periodic smooth $\mathbb{R}$-families $\{\varphi_t\}_{t\in \mathbb{R}}$ of contactomorphisms.} such that
\[ (h+\varepsilon)\# f\geqslant \ell\]
for all $\varepsilon \geqslant 0$.
\end{Lemma}
\begin{proof}
Without loss of generality, assume $\ell>0$. By the assumptions in Theorem~\ref{thm:complement-big}, {the set $A$ is the preimage of a closed subset under a smooth Reeb-invariant map $\partial W\to\mathbb{R}$. Since any closed subset of $\mathbb{R}$ can be seen as the preimage of $\{0\}$ under some smooth function $\mathbb{R}\to\mathbb{R}$,} there exists a smooth Reeb-invariant function\footnote{{Obtained as the composition.}} $\mu\colon \partial W\to\mathbb{R}$ such that $A=\mu^{-1}(\{0\})$. Since $A$ is compact and immaterial, there exists (see \cite[Lemma~7.2]{uljarevic2023selective}) a contractible loop of contactomorphisms \smash{$\varphi^f_t\colon \partial W\to\partial W$} such that
\begin{equation} \label{eq:main2ell}
\min_{x\in A, t\in\mathbb{R}} f_t(x)\geqslant 2\ell.
\end{equation}
The set
\[S:=\bigl\{ x\in\partial W\mid \bigl(\min_{t\in\mathbb{R}} f_t\bigr)(x)\leqslant \ell \bigr\}\]
is closed because the function
\[ \min_{t\in\mathbb{R}} f_t = \min_{t\in[0,1]} f_t \colon \ \partial W\to \mathbb{R} \]
is {continuous}. In addition, $S$ and $A$ are disjoint due to \eqref{eq:main2ell}. Hence, there exists $a>0$ such that $\mu^{-1}(-a, a)$ is disjoint from $S$. Let $m:=\min f$ be the (global) minimum of $f$ and let $\beta\colon \mathbb{R}\to[0,+\infty)$ be a smooth function such that
\begin{enumerate}\itemsep=0pt
\item[(1)] there exists $\delta>0$ such that $\beta(s)=0$ if, and only if, $s\in[-\delta, \delta]$, and
\item[(2)] $\beta(s)\geqslant \ell - m$ for $s\not\in (-a, a)$.
\end{enumerate}
Consequently, $\beta\circ \mu(x)\geqslant \ell-m$ for $x$ such that $\min_{t\in\mathbb{R}} f_t(x)\leqslant \ell$ and ${\rm supp} \beta\circ \mu \subset \Omega= \partial W\setminus A$. Denote $h:=\beta\circ\mu$. By the construction of $h$, we have $h+f_t\geqslant \ell$. Therefore,
\begin{align*}
((h + \varepsilon)\# f)_t(x) &= h(x) +\varepsilon + f_t\circ\bigl(\varphi^{h+\varepsilon}_t\bigr)^{-1}(x) = (h +\varepsilon + f_t)\circ\bigl(\varphi^{h+\varepsilon}_t\bigr)^{-1}(x) \geqslant \ell.
\end{align*}
Here, we used \smash{$ (h +\varepsilon )\circ\bigl(\varphi^{h+\varepsilon}_t\bigr)^{-1}= h + \varepsilon$}, which holds because $h+\varepsilon$ is time-independent and Reeb-invariant. By perturbing $h$ slightly {(within the class of non-negative smooth functions having the same support as $h$)} if necessary, one can achieve $h\in\mathcal{H}_\Omega$. This completes the proof.
\end{proof}

A contact Hamiltonian $h$ is called \emph{admissible}\label{def:admissible} if its time-1 map $\varphi^h_1$ has no fixed points $x$ such that \smash{$\bigl(\bigl(\varphi^h_1\bigr)^\ast\alpha\bigr)\bigr|_x= \alpha |_x$}. The admissibility condition is required for $HF_\ast(h)$ to be well defined. In the special case of a constant contact Hamiltonian $h\equiv a\in \mathbb{R}$, the admissibility simply means that $a$ is not a period of a Reeb orbit.
The following lemma asserts that the images of $HF_\ast(h)$ and $HF_\ast(g)$ under the canonical morphisms coincide provided $h$ and $g$ give rise to the same element in \smash{$\widetilde{{\rm Cont}}_0(\partial W)$}.

\begin{Lemma}\label{lem:compwithcontractible}
Let $h_t\colon \partial W\to \mathbb{R}$ be an admissible contact Hamiltonian and let $\varphi^f_t\colon \partial W\to\partial W$ be a contractible loop of contactomorphisms. Then, the images of the following natural morphisms are equal:
\begin{align*}
& HF_\ast(h)\to SH_\ast(W),\\
& HF_\ast(h\#f)\to SH_\ast(W).
\end{align*}
\end{Lemma}
\begin{proof}
Since $f$ gives rise to a contractible loop of contactomorphisms, there exists a {\emph{zig-zag isomorphism}} $HF_\ast(h)\to HF_\ast(h\# f)$ (see \cite[Section~9]{uljarevic2023selective} or \cite[Section~2.1.3]{djordjevic2025quantitative}). Proposition~2.4 from \cite{djordjevic2025quantitative} asserts that the following diagram, consisting of the canonical morphisms and a zig-zag isomorphism, commutes:
\[\begin{tikzcd}
&SH_\ast(W)&\\
HF_\ast(h)\arrow["\text{zig-zag}"]{rr}\arrow{ru}&&HF_\ast(h\#f).\arrow{lu}
\end{tikzcd}\]
From this commutative diagram, the claim of the lemma follows directly.
\end{proof}

\begin{proof}[Proof of Theorem~\ref{thm:complement-big}]
For an admissible contact Hamiltonian $h_t\colon \partial W\to\mathbb{R}$, denote by $\mathfrak{Im}(h)$ the image of the natural morphism $HF_\ast(h)\to SH_\ast(W)$. Let $\theta\in SH_\ast(W)$. Then, there exists a~real number $\ell\in \mathbb{R}$ that is not a period of a Reeb orbit on $\partial W$ such that $\theta\in\mathfrak{Im}(\ell)$. By Lemma~\ref{lem:twisttrick}, there exist $h\in\mathcal{H}_\Omega(\partial W)$ and a contractible loop $\varphi^f_t\colon \partial W\to\partial W$ of contactomorphisms such that the natural morphism $HF_\ast(\ell)\to SH_\ast(W)$ factors through $HF_\ast((h+\varepsilon)\#f)\to SH_\ast(W)$ for all positive $\varepsilon \in \Pi(h)$. Hence, $\mathfrak{Im}(\ell)\subset \mathfrak{Im}((h+\varepsilon)\# f)$ for all positive $\varepsilon \in \Pi(h)$. Since $f$ generates a contractible loop of contactomorphisms, Lemma~\ref{lem:compwithcontractible} implies $\mathfrak{Im}(\ell)\subset \mathfrak{Im}(h+\varepsilon)$ for all positive $\varepsilon \in \Pi(h)$. This further implies that $\theta$ belongs to the image of the continuation map $SH_\ast^\Omega(W)\to SH_\ast(W)$. Since $\theta$ was arbitrarily chosen, this completes the proof.
\end{proof}

\section{Proof of Theorem~\ref{thm:tau}}\label{sec:tau}

The quasi-measure $\tau$ is constructed in \cite{djordjevic2025quantitative} using contact spectral invariants $c(h, \theta)$ introduced independently and simultaneously in \cite{cant2023shelukhin} and \cite{djordjevic2023quantitative}.\footnote{In the present paper, we follow the approach and the conventions from \cite{djordjevic2025quantitative}} The contact spectral invariant $c(h, \theta)$ is assigned to an element $\theta\in SH_\ast(W)$ and a time-dependent contact Hamiltonian $h_t\colon \partial W\to \mathbb{R}$. {For} the construction of $\tau$, the case where $\theta$ is equal to the unit $e\in SH_\ast(W)$ is relevant.
Now, we recall the properties of the contact spectral invariants that are needed in the present paper:
\begin{enumerate}\itemsep=0pt
\item[(1)] $c(s, e)= s$ for all $s\in \mathbb{R}$.
\item[(2)] ({Invariance}) If $\bigl\{\varphi^h_t\bigr\}$ and $ \{\varphi_t^g \}$ represent the same class in \smash{$\widetilde{{\rm Cont}}_0(\partial W)$}, then $c(h,e)= c(g, e)$.
\item[(3)] (Monotonicity) If $h\leqslant g$ pointwise, then ${c(h, e)\leqslant c(g, e)}$.
\item[(4)] (Triangle inequality) If $h$ and $g$ are Reeb invariant, then\footnote{Notice that, due to our conventions, the triangle inequality is ``reversed''.}
\[
c(h, e) + c(g, e)\leqslant c(h\#g, e).
\]
\item[(5)] (Stability) If $h$ and $g$ are Reeb-invariant, then
\[ \min (h-g)\leqslant c(h,e) - c(g, e)\leqslant \max(h-g). \]
\end{enumerate}

For a Reeb-invariant contact Hamiltonian $h\colon \partial W\to\mathbb{R}$, the contact quasi-state $\zeta(h)$ is defined~by
\[\zeta(h):= \lim_{k\to+\infty} \frac{c(k\cdot h, e)}{k}.\]
The limit exists and is finite \cite[Lemma~7.2]{djordjevic2025quantitative}. The contact quasi-state $\zeta$ satisfies the following homogeneity property: $\zeta(sh)= s \zeta(h)$ for $s\geqslant 0$. Now, the contact quasi-measure $\tau$ is defined by
\[ \tau(A):= \inf \{ \zeta(h)\mid h\colon\partial W\to[0,1] \text{ Reeb-invariant, smooth, and } h|_A=1 \} \]
for a closed subset $A\subset\partial W$. The following lemma is used in the proof of Theorem~\ref{thm:tau}.

{
\begin{Lemma}\label{lem:aut}
Let $M$ be a closed contact manifold with a~contact form $\alpha$ and let $B\subset M$ be a~closed subset. Let $f_t\colon M\to\mathbb{R}$ be a~contact Hamiltonian and $h\colon M\to[0,1]$ a Reeb-invariant contact Hamiltonian such that $ h |_B=1$. Denote
\[\ell:=\inf_{t\in \mathbb{R}, x\in M\setminus B} f_t(x)\qquad \text{and} \qquad
m:=\inf_{t\in \mathbb{R}, x\in M} f_t(x),\]
and assume $m$ is finite.\footnote{This implies $\ell$ is finite as well.} Then, the following inequality holds:
\[ ((\ell-m)\cdot h)\# f\geqslant \ell. \]
\end{Lemma}}
\begin{proof}
Denote $g:=(\ell-m)\cdot h$. Since $g$ is time-independent and Reeb-invariant, it is constant along its trajectories. Consequently,
\begin{align*}
(g\# f)_t(x) = g(x) + f_t\circ (\varphi^{g}_t )^{-1}(x) = (g + f_t)\circ (\varphi^{g}_t )^{-1}(x).
\end{align*}
Therefore, it is enough to prove $g+f_t\geqslant \ell$. If $x\in {M\setminus B}$, then $f_t(x)\geqslant \ell$. Since $g$ is non-negative, this implies $g(x)+ f_t(x)\geqslant \ell$. If, on the other hand, $x\not \in {M\setminus B}$, then $x\in B$ and
\[ g(x) + f_t(x)= \ell-m + f_t(x)\geqslant \ell. \]
This finishes the proof.
\end{proof}

\begin{proof}[Proof of Theorem~\ref{thm:tau}]
Let $h\colon\partial W\to[0,1]$ be a Reeb-invariant contact Hamiltonian such that~${ h |_B=1}$. By Lemma~\ref{lem:aut}, we get $ ((m_B-m)\cdot h )\#f\geqslant m_B$. The {invariance} and the monotonicity properties of the contact spectral invariants imply
\[c((m_B-m)\cdot h, e)= c ( ((m_B-m)\cdot h )\#f, e )\geqslant c(m_B, e)=m_B.\]
The same argument applied to the iterated loop $\varphi^f_{kt}$ (generated by the contact Hamiltonian ${k\cdot f_{kt}}$), $k\in\mathbb{N}$, implies
\[c(k\cdot(m_B-m)\cdot h, e)\geqslant k\cdot m_B.\]
Therefore,
\[\zeta((m_B-m)\cdot h)= \lim_{k\to+\infty} \frac{c(k\cdot(m_B-m)\cdot h, e)}{k}\geqslant m_B. \]
Consequently, due to homogeneity of $\zeta$, we have \smash{$\zeta(h)\geqslant\frac{m_B}{m_B-m}$}. This completes the proof.
\end{proof}

\subsection*{Acknowledgements}

The results of this paper rely on selective symplectic homology and contact quasi-measures. Both theories were inspired to a great extent by joint work with Will Merry \cite{merry2019maximum}. I am grateful to Will Merry for his friendship and for this collaboration. I would like to thank the referees for their thorough review. The author is partially supported by the Ministry of Science, Technological Development and Innovation, Republic of Serbia, through the project 451-03-33/2026-03/200104.


\pdfbookmark[1]{References}{ref}
\LastPageEnding

\end{document}